\documentclass[a4paper, 12pt]{article}

\usepackage[koi8-r,cp1251]{inputenc}
\usepackage{amsfonts,amssymb,amsmath, hyperref}
\usepackage[final]{epsfig}
\usepackage{graphicx}

\begin{document}
\begin{center}
\Large{On  contractibility of the Gelfand spectrum of semigroup measure algebras}
\end{center}

\begin{center}
A. R. Mirotin
\end{center}

\begin{center}
amirotin@yandex.ru
\end{center}

\

Abstract.  \small{Sufficient conditions  for a semigroup measure algebra to have contractible  Gelfand spectrum are given and it is shown that for a wide class of semigroups  these conditions are also necessary.}

\


\

In the paper \cite{S} by A. Sasane it was shown that the Gelfand spectrum (equipped with the Gelfand topology) of the  algebra of complex bounded  measures on the additive semigroup $[0,+\infty)$  is contractible. Several applications of this result are mentioned  in \cite{S}, too. The main goal of this note is to give a  simple proof of a generalization of the main result of \cite{S} to general semigroups. Actually, we give sufficient conditions  for a semigroup measure algebra $\mathcal{M}(S)$  to have contractible  Gelfand spectrum $\Delta\mathcal{M}(S)$  and show that for a wide class of semigroups  these conditions are also necessary. On the other hand, we show that for a nontrivial  locally compact Abelian group $G$  the space  $\Delta\mathcal{M}(G)$   is noncontractible.

  In the following $S$ stands for an Abelian semitopological locally compact semigroup with unit $e,$ and let $\mathcal{M}(S)$ be the algebra of all complex bounded regular Borel measures on $S$ with  respect to the  convolution  (see, e.g. \cite[Chapter 5]{GMcG})
$$
\int_Sfd\mu\ast\nu:=\int_{S\times S}f(xy)d\mu(x)d\nu(y)\ (f\in C^b(S)).
$$

  Let $S^*_1$ denotes the semigroup  of  nonnegative bounded  semicharacters of $S$ (nontrivial homomorphisms to the multiplicative semigroup $[0,1]$) that are $\mu$- measurable for all $\mu\in \mathcal{M}(S)$ endowed with the topology of pointwise convergence.

\textbf{Theorem 1.}  Let $S$ be an Abelian  locally compact semitopological semigroup with unit $e.$

1) If $S$ has no nontrivial invertible elements and the space $S^*_1$ is  path-connected then the Gelfand spectrum $\Delta\mathcal{M}(S)$  of the algebra $\mathcal{M}(S)$  is contractible.

2) Let  $S$ has the form $P\cup\{e\}$ where $P$ is an open subsemigroup of a locally compact  Abelian   group $G$ and the unit $e$ of $G$ belongs to the closure of $P.$ If  $\Delta\mathcal{M}(S)$ is contractible, the  space $S^*_1$ is  path-connected.

\textit{Proof.} 1)  (cf. \cite{S}). Denote by $1_e$ the characteristic function of the set $\{e\}.$ Clearly, $1_e \in S^*_1.$ Let $(\rho_t)_{t\in [0,1]}$ be a continuous  path in $S^*_1$ such that $\rho_0=\textbf{1}, \rho_1=1_e$  ($\textbf{1}$ denotes the unit semicharacter of $S.$) Define a continuous map $H:\Delta\mathcal{M}(S)\times [0,1]\to \Delta\mathcal{M}(S)$ as follows:
 $$
H(\varphi,t)(\mu)=\varphi(\rho_t\cdot\mu)
$$
 where $\rho_t\cdot\mu$ means the product of a function and a measure  in the sense of \cite{Burb}. We show that $H$ is well-defined. Indeed, since $\rho_t\cdot(\mu\ast\nu)=(\rho_t\cdot\mu)\ast(\rho_t\cdot\nu)$
(see, e.g., \cite[Chapter VIII, Section 3,  Proposition 6]{Burb}), we have $H(\varphi,t)\in \Delta\mathcal{M}(S).$ Moreover, by the inequality ($\varphi, \varphi_0\in \Delta\mathcal{M}(S),$  $\mu\in \mathcal{M}(S)$)
$$
|H(\varphi,t)(\mu)-H(\varphi_0,t_0)(\mu)|\leq |\varphi(\rho_t\cdot\mu)-\varphi(\rho_{t_0}\cdot\mu)|+|\varphi(\rho_{t_0}\cdot\mu)-\varphi_0(\rho_{t_0}\cdot\mu)|
$$
 the map $H$ is continuous because
$$
|\varphi(\rho_t\cdot\mu)-\varphi(\rho_{t_0}\cdot\mu)|\leq \|\varphi\|\|\rho_t\cdot\mu-\rho_{t_0}\cdot\mu\|=\|(\rho_t-\rho_{t_0})\cdot\mu\|\leq \int_S |\rho_t-\rho_{t_0}|d|\mu|\to 0
$$
as $t\to t_0$ by the Lebesgue theorem.

Finally, if we define $\varphi_1(\mu):=\mu(\{e\}),$ then  $H(\varphi,0)=\varphi,$ and $H(\varphi,1)=\varphi_1,$ for all $\varphi\in  \Delta\mathcal{M}(S).$ So, $\Delta\mathcal{M}(S)$
is contractible.

2)     Since  $\Delta\mathcal{M}(S)$ is contractible,  there exists a  continuous map $H:\Delta\mathcal{M}(S)\times [0,1]\to \Delta\mathcal{M}(S)$ such that $H(\varphi,0)=\varphi,$ and $H(\varphi,1)=\varphi_1$ for all $\varphi\in  \Delta\mathcal{M}(S)$ where $\varphi_1(\mu):=\mu(\{e\})$. Let $\varphi_0(\mu):=\int_Sd\mu\ (\varphi_0\in \Delta\mathcal{M}(S)).$ For $t\in [0,1]$ define $\rho_t:S\to \mathbb{R}_+$ as follows:
$$
\rho_t(s)=|H(\varphi_0,t)(\delta_s)|
$$
where  $\delta_s$ denotes the Dirac measure centered on $s\in S.$ Then  $\rho_t(e)=1$ and $\rho_t(s)\leq\|H(\varphi_0,t)\|\|\delta_s\|=1.$ Since $\delta_a\ast\delta_b=\delta_{ab}\ (a,b\in S),$ the map $\rho_t$ is a homomorphism from $S$ to $[0,1].$
     By the result of Devinatz and Nussbaum \cite{DN} the restriction $\rho_t|P$ is continuous and we conclude that $\rho_t\in S^*_1$. Since the map $t\mapsto H(\varphi_0,t)$ is continuous, we get a continuous path $(\rho_t)_{t\in[0,1]}$ in $S^*_1$ such that $\rho_0=\textbf{1}, \rho_1=1_e$. Therefore for every $\rho\in S^*_1$ we have a continuous path $(\rho\rho_t)_{t\in[0,1]}$ in $S^*_1$ from $\rho$ to $1_e$ which completes the proof.

\textbf{Remark 1.}  The condition that $S$ has no nontrivial invertible elements is essential in theorem 1 as  examples  of locally compact Abelian groups show (see  theorem 2 below). We conjecture that this condition is in fact necessary for $\Delta\mathcal{M}(S)$  to be contractible.

\textbf{Remark 2 }(cf. \cite{S}).  The part 1 of  theorem 1 remains valid (along with its proof) for every  unital Banach subalgebra $\mathcal{R}$ of $\mathcal{M}(S)$ with the
property: for some continuous path $(\rho_t)_{t\in[0,1]}$ in $S^*_1$ from $\textbf{1}$ to $1_e$
 $$
 \forall\mu\in \mathcal{R},\forall t\in[0,1]\quad \rho_t\cdot\mu\in \mathcal{R}.\eqno(1)
 $$

\textbf{Examples.} Let $S$ be an Abelian topological semigroup with invariant measure $m$ \cite{Mir}  and unit $e$  without nontrivial invertible elements (in particular, let $S$ be a Borel subsemigroup of a locally compact Abelian group with nonempty interior and without nontrivial subgroups, $e\in S$). Then the property (1) holds for the subalgebra  $\mathcal{M}_d(S)$ of atomic measures on $S$, and for the  subalgebra  $L^1(S,m)+\mathcal{M}_d(S)$  of all complex Borel measures that do not have a singular non-atomic part with respect to $m$; the property (1)  holds also for the  algebra $L^1(S,m)+\mathbb{C}\delta_e.$ So the Gelfand spectra of these algebras  are contractible. The same is true for images of aforementioned algebras with respect to the generalized Laplace transform in the sense of N. Bourbaki if the semigroup $\widehat{S}$ of all continuous bounded semicharacters of $S$ is complete  (see \cite[Chapter IX, subsection 5.7]{Burb}, and also  \cite{UMZ1992}, \cite{VUZ1995}).

\textbf{Corollary 1} (cf. \cite{S}).  Let $S$ has no nontrivial invertible elements and the space $S^*_1$ is  path-connected. Every unital Banach subalgebra $\mathcal{R}$ of $\mathcal{M}(S)$ with the property (1)
  is a Hermite ring in a sense that for every  left invertible $n\times k$ matrix $(f_{ij})$ with  entries in  $\mathcal{R}$ ($k<n$) there exists an invertible $n\times n$ matrix $(g_{ij})$ with  entries in  $\mathcal{R}$ such that $g_{ij}=f_{ij}$ for  all $1\leq i \leq n$ and $1\leq j \leq k$.

It follows from the remark 1 and  \cite[Theorem 3, p. 127]{Lin}.

\textbf{Remark 3.} For  applications of results like corollary 1 to  the problems in control theory see \cite{S} and references there.

\textbf{Corollary 2.} Let $S$ be an Abelian  locally compact semitopological semigroup with unit $e$ and without nontrivial invertible elements. If there is  $\rho\in S^*_1$ such that $0<\rho(s)<1$ for all $s\in S\setminus\{e\}$ then for every unital Banach subalgebra $\mathcal{R}$ of  $\mathcal{M}(S)$ with the property (1) the Gelfand spectrum $\Delta\mathcal{R}$  is contractible.

Indeed, in this case we have a  continuous path  $(\rho_t)_{t\in[0,1]}$ in $S^*_1,$  where $\rho_t:=\rho^{t/(1-t)}$ for $t\in[0,1),$ $\rho_1:=1_e,$ from $\textbf{1}$ to $1_e$.

\textbf{Theorem 2.} For a nontrivial  locally compact Abelian group $G$  the space  $\Delta\mathcal{M}(G)$   is noncontractible.

Proof. Suppose that the Gelfand spectrum  $\Delta=\Delta\mathcal{M}(G)$ is contractible.  Then the one-dimensional \v{C}ech cohomology group $H^1(\Delta)$ with integral coefficients is trivial (see, e.g., \cite[Chapter 2, section 2.1]{Massey}).  On the other hand, since $L^1(G_d)$ is a maximal group algebra in $\mathcal{M}(G)$ \cite[8.2.1]{Tay} ($G_d$ denotes the group $G$ endowed with discrete topology), Theorem 8.1.3 from \cite{Tay} (see also \cite[Theorem 5.3.7]{GMcG}) implies that $H^1(\Delta)$ contains the isomorphic image of
$H^1(\Delta L^1(G_d))=H^1(\widehat{G_d}).$ If the dual group $\widehat{G_d}$ is connected, it is known \cite[8.3.2]{Tay} that $H^1(\widehat{G_d})$ is isomorphic to $\widehat{\widehat{G_d}}=G$ and we have a contradiction with $H^1(\Delta)=0$. Now let $\widehat{G_d}$ be disconnected. Then the group $G_d$ contains a nontrivial finite subgroup $F$ (see, e.g., \cite[(24.19)]{HiR}). If $F\ne G$ the normalized Haar measure of $F$ is a nontrivial idempotent in  $\mathcal{M}(G)$ and therefore $\Delta$ is disconnected by Shilov's Idempotent Theorem  (and if  $G=F$ this is obvious). This completes the proof.

We say that $S$ is an \textit{$M$-semigroup} if every its semicharacter  of the form $s\mapsto |\varphi(\delta_s)|\ (\varphi\in \mathcal{M}(S))$ belongs to $S^*_1.$ The proof of the part 2 of theorem 1 and results of \cite[section 3]{Baartz} give us examples of $M$-semigroups.  On the other hand,  there are locally compact Abelian topological semigroups of idempotents which are not $M$-semigroups \cite[subsection 2.4]{Baartz}.

\textbf{Theorem 3.} Let $S$ be an Abelian  locally compact semitopological $M$-semigroup with unit $e$   and without nontrivial invertible elements. Then $\Delta\mathcal{M}(S)$    is contractible if  and only if  the space $S^*_1$ is  path-connected.

Indeed, the sufficiency  follows from part 1 of theorem 1. As regards  the necessity, it can be    proved in the same way as  the part 2 of theorem 1.

\textbf{Acknowledgments.}
 The paper  was published in Semigroup Forum - 2019 - Vol. 98, No 1 - P. 209 -212;
  https://doi.org/10.1007/s00233-018-9974-x.
The author is deeply indebted to the referee for very useful comments and suggestions that improve the paper.

\end{document}